\documentclass[12pt]{article}

\usepackage{amsthm,amssymb}
\usepackage{amsmath,amsfonts,latexsym}
\usepackage{color}
\usepackage{float}

\newtheorem{defi}{Definition}

\newtheorem{tm}{Theorem}

\newtheorem{kor}{Corollary}

\newtheorem{rem}{Remark}

\begin{document}

\begin{titlepage}

\vspace*{0.2cm}

\begin{center}
{\Large\bf  Self-orthogonal codes from equitable partitions of association schemes}       
\end{center}

\vspace*{0.2cm}

\begin{center}
Dean Crnkovi\'c  \\
({\it\small E-mail: deanc@math.uniri.hr})\\[3pt] 
Sanja Rukavina \\
({\it\small E-mail: sanjar@math.uniri.hr})\\[3pt]
and\\[3pt]
Andrea \v Svob \\
({\it\small E-mail: asvob@math.uniri.hr})\\[3pt]
{\it\small Department of Mathematics, University of Rijeka} \\
{\it\small Radmile Matej\v ci\'c 2, 51000 Rijeka, Croatia}\\
\end{center}

\vspace*{0.5cm}

\begin{abstract}
We give a method of constructing self-orthogonal codes from equitable partitions of association schemes.
By applying this method we construct self-orthogonal codes from some distance-regular graphs.
Some of the obtained codes are optimal.
Further, we introduce a notion of self-orthogonal subspace codes. 
We show that under some conditions equitable partitions of association schemes yield such self-orthogonal subspace codes and we give some examples from distance-regular graphs.
\end{abstract}

\vspace*{0.5cm}



\vspace*{1.5cm}

\begin{center}{\bf Acknowledgement}\end{center}
This work has been fully supported by {\rm C}roatian Science Foundation under the project 6732.

\end{titlepage}

\section{Introduction}\label{intro_new}

In this paper we give a construction of self-orthogonal codes from equitable partitions of association schemes. 
The results presented in the paper can be seen as a generalization of results on construction of self-orthogonal codes from orbit matrices of 2-designs (see \cite{2-des-codes, har_ton}) 
and strongly regular graphs (see \cite{SRG-codes}). Moreover, we study subspace codes obtained from equitable partitions of association schemes and introduce a notion of self-orthogonal subspace codes.

Self-orthogonal codes have wide applications in communications (see \cite{RainsSloane}), including, for example, in secret sharing. 
Further, some of the most interesting and best linear codes known are self-orthogonal, \textit{e.g.} the $[8,4,4]$ Hamming code, the $[24,12,8]$ extended binary Golay code
and the $[12,6,6]$ ternary Golay code. These are the reasons why our interest is also in defining and constructing self-orthogonal subspace codes. 
Subspace codes are relatively new topic of interest (see \cite{network-book}) giving a new approach to network coding. 

Beside the theoretical results on a construction of self-orthogonal codes and self-orthogonal subspace codes from equitable partitions of association schemes,
we also present examples obtained by applying the described methods on association schemes related to some distance-regular graphs. Some of the obtained linear codes are optimal.

The paper is outlined as follows. In the next section we provide the relevant background information. In Section \ref{SO_linear} we give a construction method. 
Following this, we describe our construction using equitable partitions of distance-regular graphs. In Section \ref{SO-subspace} we introduce a notion of self-orthogonal subspace codes, 
giving a method for obtaining such subspace codes. We conclude with our construction of self-orthogonal subspace codes via distance-regular graphs and give examples to demonstrate the construction.

In this work we have used computer algebra systems GAP \cite{GAP2016}, Magma \cite{magma} and Sage \cite{sage}, and Hanaki's programs \cite{Hanaki}.

\section{Preliminaries} \label{intro}

We assume that the reader is familiar with the basic facts of theory of distance-regular graphs and association schemes. 
For background reading in theory of distance-regular graphs and association schemes we refer the reader to \cite{BCN}. For further reading on the topic we refer the reader to \cite{GoM} and \cite{BMTh}.
We also assume a basic knowledge of coding theory (see \cite{FEC}).\\

We will follow the definition of an association scheme given in \cite{BCN}, although some authors use a term a {\it symmetric association scheme} for such structure.

Let $X$ be a finite set. An {\it association scheme} with $d$ classes is a pair $(X,\mathcal{R})$ such that
\begin{enumerate}
 \item $\mathcal{R}= \{ R_0,R_1,...,R_d\}$ is a partition of $X \times X$,
 \item $R_0=\bigtriangleup = \{ (x,x) |x \in X \} $,
 \item $R_i=R_i^T$ (\textit{i.e.} $(x,y) \in R_i \Rightarrow (y,x) \in R_i)$ for all $i \in \{ 0,1,...d \} $,
 \item there are numbers $p_{ij}^k$ (the intersection numbers of the scheme) such that for any pair $(x,y) \in R_k$ the number of $z \in X$ such that $(x,z) \in R_i$ and $(z,y) \in R_j$ equals $p_{ij}^k$. 
\end{enumerate}

The relations $R_i$, $i \in \{ 0,1, ..., d\}$, of an association scheme can be described by the set of symmetric $(0,1)$-adjacency matrices $\mathcal{A}= \{A_0, A_1, ..., A_d \}$, $[A_i]_{xy}=1$ if 
$(x,y) \in R_i$, which generate $(d+1)$-dimensional commutative and associative algebra over real or complex numbers called {\it the Bose-Mesner algebra} of the scheme. 
The matrices $\{A_0, A_1, ..., A_d \}$ satisfy
\begin{equation} \label{form1}
 A_i A_j=\sum_{k=0}^d p_{i,j}^k A_k=A_j A_i.
\end{equation}
Each of the matrices $A_i$, $i \in \{1,2, ..., d\}$, represents a simple graph $\Gamma_i$ on the set of vertices $X$ (if $(x,y)\in R_i$ then vertices $x$ and $y$ are adjacent in $\Gamma_i$), 
and the graphs $\Gamma_i$ form an edge-coloring of the complete graph on $X$.

A $q$-ary \emph{linear code} $C$ of dimension $k$ for a prime power $q$, is a $k$-dimensional subspace of a vector space $\mathbb{F}_{q}^n$.  Elements of $C$ are called codewords.
Let $x=(x_1,...,x_n)$ and $y=(y_1,...,y_n)\in \mathbb{F}_q^n$. The \emph{Hamming distance} between words $x$ and $y$ is the number \ $d(x,y)=\left| \{ i  :  x_i \neq y_i \} \right| $.
The \emph{minimum distance} of the code  $C$ is defined by \ $d=\mbox{min}\{ d(x,y):  x,y\in C, \ x\neq y\}$.  The \emph{weight} of a codeword $x$ is \ $w(x)=d(x,0)=|\{i  :  x_i\neq 0 \}|$.  
For a linear code, $d=\mbox{min}\{ w(x)  :  x \in C, x\neq0 \}.$

The \emph{dual} code $C^\perp$ of a code $C$ is the orthogonal complement of $C$ under the standard inner product
$\langle\cdot \, ,\cdot\rangle$, i.e.\ $C^{\perp} = \{ v \in \mathbb{F}_{q}^n| \langle v,c\rangle=0 {\rm \ for\ all\ } c \in C \}$. A code $C$ is \emph{self-orthogonal} if $C \subseteq C^\perp$ and 
\emph{self-dual} if $C = C^\perp$.

A linear code is called \textit{projective} if no two columns of the generator matrix are linearly dependent. 
A code is projective if and only if the minimum distance of its dual code is at least three.
A \textit{two-weight code} is a code which has only two non-zero weights. 
A projective two-weight codes are related to strongly regular graphs.

A  $q$-ary linear code of length $n$, dimension $k$, and distance $d$ is called a $[n,k,d]_q$ code. 
An $[n, k]$ linear code $C$ is \emph{optimal} if the minimum weight of $C$ achieves the theoretical upper bound on the minimum weight of $[n, k]$ linear codes, and \emph{near-optimal} if its minimum weight is at most 1 less than the largest possible value.  An $[n, k]$ linear code $C$ is said to be a \emph{best known} linear $[n, k]$ code if $C$ has the highest minimum weight among all known 
$[n, k]$ linear codes.  A catalogue of best known codes is maintained at \cite{codetables}, to which we compare the minimum weight of all codes constructed in this paper.

Linear codes over finite rings are defined similarly to finite fields, where the codes are modules instead of vector spaces. 
$\mathbb{Z}_m $ denotes the ring of integers modulo $m$, where $m$ is a positive integer, $m\geq 2$.
The most notable codes over rings are codes over $\mathbb{Z}_4$.

\bigskip

A {\it subspace code} $C_S$ is a nonempty set of subspaces of $\mathbb F^n_q$. 
R. K\"{o}tter and F. Kschischang proved (see \cite{network-coding}) that subspace codes are efficient for transmission in networks.
For the parameters of a subspace code we will follow the notation from \cite{table_sc}
and use a {\it subspace distance} given by
\begin{equation} \label{submet}
d_s(U,W)=dim(U+W)-dim(U\cap W),
\end{equation}
where $U, W \in C_S$. The {\it minimum distance} of $C_S$ is given by 
$$d=min \{d_S(U,W) | \ U,W \in C_S, U \neq W \}.$$
A code $C_S$ is called an $(n,\#C_S,d;K)_q$ subspace code if the dimensions of the codewords of $C_S$ are contained in a set $K\subseteq \{0,1,2, ..., n\}$. 
In the case $K=\{k\}$, a subspace code $C_S$ is called a {\it constant dimension code} with the parameters $(n,\#C_S,d;k)_q$, otherwise, \textit{i.e.} if all codewords do not have the same dimension, 
$C_S$ is called a {\it mixed dimension code}. Such subspace code is denoted by $(n,\#C_S,d)_q$.
For reading on recent results on subspace codes we refer the reader to \cite{network-book, heinlein-dcc}. 

In this paper we give a construction of self-orthogonal codes using equitable partitions of association schemes.
In the case of 2-class association schemes this method of construction coincides with the construction from strongly regular graphs given in \cite{SRG-codes}.
We apply this method to construct self-orthogonal codes from some distance-regular graphs of diameter $d$, $3\leq d\leq 9$, \textit{i.e.}  Doro graph of diameter $d=3$, 
Hadamard graph on 48 vertices of diameter $d=4$, Doubled Gewirtz graph of diameter $d=5$, Incidence graph of $GH(3,3)$ of diameter $d=6$, 
Doubled Odd graph $D(O_4)$ of diameter $d=7$ and Foster graph of diameter $d=8$. Some of the obtained self-orthogonal codes are optimal \textit{i.e.} they reach the theoretical upper bound. 
Further, for obtaining self-orthogonal subspace codes we apply the method on distance-regular graph called Doubled Higman-Sims of diameter $d=5$.

\section{Self-orthogonal linear codes from equitable partitions of association schemes} \label{SO_linear}

Suppose $A$ is a symmetric real matrix whose rows and columns are indexed by the elements of $X= \{1, \ldots ,n \}$. Let $\{C_0,...,C_{t-1} \}$ be a partition of $X$. 
The characteristic matrix $H$ is 
the $n \times t$ matrix whose $j$th column is the characteristic vector of $C_j$, where $j= 0,...,t-1$. 

A partition $\Pi =\{ C_0,C_1,...,C_{t-1} \}$ of the $n$ vertices of a graph $G$ is {\it equitable} (or regular) if for every pair of (not necessarily distinct) indices
$i,j \in \{0,1,...,t-1 \}$ there  is  a  nonnegative  integer $b_{i,j}$ such  that  each  vertex $v \in C_i$ has  exactly $b_{i,j}$ neighbors in $C_j$, regardless of the choice of $v$.  
The $t \times t$ quotient matrix $B = (b_{i,j})$ is well-defined if and only if the partition $\Pi$ is equitable.
An equitable (or regular) partition of an association scheme $(X,\mathcal{R})$ is a partition of $X$ which is equitable with respect to each of the graphs $\Gamma_i$, $i \in \{1,2, ...,d \}$
corresponding to the association scheme $(X,\mathcal{R})$ with $d$ classes.

Let $\Pi$ be an equitable partition of $(X,\mathcal{R})$ with $t$ cells, and let $H$ be the characteristic matrix of the partition $\Pi$.
Further, let $A_i$ be the adjacency matrix corresponding to a relation $R_i$.
Then the following holds:
\begin{equation} \label{form2}
 A_iH=HM_i,
\end{equation}
where $M_i$ denotes the corresponding $t \times t$ quotient matrix of $A_i$ with respect to $\Pi$.
The matrix $H^T H$ is diagonal and invertible and, therefore,
\begin{equation} \label{form2a}
M_i=(H^T H)^{-1} H^T A_iH.
\end{equation}

\begin{tm} \label{th_struc_const}
Let $\Pi$ be an equitable partition of a $d$-class association scheme $(X,\mathcal{R})$ with $t$ cells, 
and let $M_i$, $i=0, 1, \dots, d$, denote the quotient matrix of the graph $\Gamma_i$ with respect to $\Pi$. Then
\begin{equation} \label{form4}  
 M_i M_j=\sum_{k=0}^d p_{i,j}^k M_k=M_j M_i,
\end{equation}
where numbers $p_{ij}^k$ are the intersection numbers of the scheme.
\end{tm}
{\bf Proof}
From the equation (\ref{form2a}) it follows that 
$HM_i=H(H^T H)^{-1} H^T A_i H,$ and by applying (\ref{form2}) we get
$H(H^T H)^{-1} H^T=I.$

Now we have
\begin{align*}
M_i M_j&=(H^T H)^{-1} H^T A_i H (H^T H)^{-1} H^T A_j H\\
&= (H^T H)^{-1} H^T A_i A_j H = (H^T H)^{-1} H^T A_j A_i H = M_j M_i,
\end{align*}

and
\begin{align*}
M_i M_j&= (H^T H)^{-1} H^T A_i A_j H= (H^T H)^{-1} H^T (\sum_{k=0}^d p_{i,j}^k A_k) H\\
&= \sum_{k=0}^d p_{i,j}^k [ (H^T H)^{-1} H^T A_k H ] =  \sum_{k=0}^d p_{i,j}^k M_k.
\end{align*}
{$\Box$}
 
\begin{tm} \label{so_code}
Let $\Pi$ be an equitable partition of a $d$-class association scheme $(X,\mathcal{R})$ with $n$ cells of the same length $\frac{|X|}{n}$ and let $p$ be a prime number.
If for fixed $i \in \{1,2,...,d\}$ and for all $k \in \{0,1,...,d\}$ the number $p$ divides $p_{i,i}^k$, 
then the rows of the matrix $M_i$ span a self-orthogonal code of length $n$ over the field $\mathbb F_q$, $q=p^m$, where $m$ is a positive integer.
\end{tm}

{\bf Proof}
Since the partition $\Pi$ is equitable  
having all cells of the same size, and the matrices $A_i$, $i=0, \ldots , d$, are symmetric, the matrices $M_i$, $i=0, \ldots , d$, are symmetric too.
Now, from (\ref{form4}) it follows that
\begin{equation} \label{form6}  
 M_iM_i^T=M_i M_i=\sum_{k=0}^d p_{i,i}^k M_k.
\end{equation}
Since $p|p_{i,i}^k$ for all $k \in \{0,1, \ldots ,d\}$, $p_{i,i}^k=0$ over $\mathbb F_q$, and the statement of the theorem holds. {$\Box$} \\

A similar result can be proven for codes over rings $\mathbb{Z}_{m}$.

\begin{tm} \label{so_code-rings}
Let $\Pi$ be a $d$-class equitable partition of an association scheme $(X,\mathcal{R})$ with $n$ cells of the same length $\frac{|X|}{n}$ and let $m$ be a positive integer.
If for fixed $i \in \{1,2,...,d\}$ and for all $k \in \{0,1,...,d\}$ the integer $m$ divides $p_{i,i}^k$, 
then rows of the matrix $M_i$ span a self-orthogonal code of length $n$ over the ring $\mathbb{Z}_{m}$.
\end{tm}

\subsection{Codes from distance-regular graphs}

Let $\Gamma$ be a graph with diameter $d$, and let $\delta(u,v)$ denote the distance between vertices $u$ and $v$ of $\Gamma$.
The $i$th-neighborhood of a vertex $v$ is the set $\Gamma_{i}(v) = \{ w : \delta(v,w) = i \}$. 
Similarly, we define $\Gamma_{i}$ to be the $i$th-distance graph of $\Gamma$, that is, the vertex set of $\Gamma_{i}$ is the same as for $\Gamma$, with adjacency in 
$\Gamma_{i}$ defined by the $i$th distance relation in $\Gamma$.
We say that $\Gamma$ is distance-regular if the distance relations of $\Gamma$ give the relations of a $d$-class association scheme, that is, for every choice of $0 \leq i,j,k \leq d$, 
all vertices $v$ and $w$ with $\delta(v,w)=k$ satisfy $|\Gamma_{i}(v) \cap \Gamma_{j}(w)| = p^{k}_{ij}$ for some constant $p^{k}_{ij}$.
In a distance-regular graph, we have that $p^{k}_{ij}=0$ whenever $i+j < k$ or $k<|i-j|$.
A distance-regular graph $\Gamma$ is necessarily regular with degree $p^{0}_{11}$; more generally, each distance graph $\Gamma_{i}$ is regular with degree $k_{i}=p^{0}_{ii}$.
An equivalent definition of distance-regular graphs is the existence of the constants $b_{i}=p^{i}_{i+1,1}$ and $c_{i}= p^{i}_{i-1,1}$ for $0 \leq i \leq d$ (notice that $b_{d}=c_{0}=0$).
The sequence $\{b_0,b_1,\dots,b_{d-1};c_1,c_2,\dots,c_d\}$, where $d$ is the diameter of $\Gamma$, is called the intersection array of $\Gamma$. 
Clearly, $b_0=k$, $b_d=c_0=0$, $c_1=0$.

Let $\Gamma$ be a distance-regular graph with diameter $d$ and adjacency matrix $A$, and let $A_i$ denotes the distance-$i$ matrix of $\Gamma$, $i=0, 1, \dots, d$. 
Further, let $G$ be an automorphism group of $\Gamma$. 
In the sequel, the quotient matrix of $A_i$, $i=1, 2, \dots, d$, with respect to the orbit partition induced by $G$ will be denoted by $M_i$.

Since a distance-regular graph induces an association scheme, Theorem \ref{so_code} implies the following corollary.
This result is a generalization of the result on construction of self-orthogonal codes from orbit matrices of strongly regular graphs given in \cite{SRG-codes}. 

\begin{kor} \label{cor_drg}
Let $\Gamma$ be a distance-regular graph with diameter $d$, and let $p$ be a prime number. 
Further, let an automorphism group $G$ acts on $\Gamma$ with $n$ orbits of the same length. 
If there exists $i \in \{1,2,...,d\}$ such that
for all $k \in \{0,1,...,d\}$ the prime $p$ divides $p_{i,i}^k$, then the rows of the matrix $M_i$ span a self-orthogonal code of length $n$ over the field $\mathbb F_q$, $q=p^m$, 
where $m$ is a positive integer.
\end{kor}

The following statement is a direct consequence of Theorem \ref{so_code-rings}.

\begin{kor} \label{cor_drg-rings}
Let $\Gamma$ be a distance-regular graph with diameter $d$, and let $m$ be a positive integer. 
Further, let an automorphism group $G$ acts on $\Gamma$ with $n$ orbits of the same length. 
If there exists $i \in \{1,2,...,d\}$ such that
for all $k \in \{0,1,...,d\}$ the integer $m$ divides $p_{i,i}^k$, then the rows of the matrix $M_i$ span a self-orthogonal code of length $n$ over the ring $\mathbb{Z}_{m}$.
\end{kor}

\begin{rem}
The trivial group acts on every distance-regular graph with all orbits of the length $1$, and in that case the corresponding quotient matrix is actually the adjacency matrix of the graph.
\end{rem}

\section{Examples of self-orthogonal codes}

In this Section we give examples of self-orthogonal codes obtained by applying the method described in Theorem \ref{so_code} and Corollary \ref{cor_drg}. 
The examples of the construction from strongly regular graphs are given in \cite{SRG-codes}.
Here we apply the method by taking distance-regular graphs of diameter $d$, $3\leq d\leq 8$, one example for each diameter. 
In each subsection we give basic information on the distance-regular graphs, the coefficients $p_{i,i}^k$ of adjacency matrices and obtained self-orthogonal codes. 
More information about the each graph can be seen in \cite{BCN}. The obtained self-orthogonal codes with $*$ are optimal. Some of the obtained codes are self-dual.

\subsection{Doro graph, $d=3$}

The Doro graph, which we will denote by $\Gamma_D$, has 68 vertices and the intersection array $\{12,10,3;1,3,8\}$. 
The distance-$i$ matrices $A_i$ of $\Gamma_D$ ($i=0, 1, \dots, 3$) determine an 3-class association scheme and form a basis of a Bose-Mesner algebra. 
The coefficients of $A_i \cdot A_i$, \textit{i.e.} intersection numbers $p_{i,i}^k$, are given in Table \ref{table-coeff-D}.

\begin{table}[H] 
\begin{center} \begin{scriptsize}
\begin{tabular}{|c|c|c|c|c|}
 \hline 
 & $A_0$ & $A_1$& $A_2$ & $A_3$ \\
 \hline \hline
$A_0 \cdot A_0$ & 1 & 0 & 0 & 0 \\
 \hline
$A_1 \cdot A_1$ & 12 & 1 & 3 & 0 \\
 \hline
$A_2 \cdot A_2$ & 40 & 20 & 24 & 24 \\
 \hline
 $A_3 \cdot A_3$ & 15 & 5 & 3 & 2  \\
 \hline \hline
\end{tabular} \end{scriptsize}
\caption{\footnotesize The coefficients $p_{i,i}^k$ of adjacency matrices for the Doro graph $\Gamma_D$} \label{table-coeff-D}
\end{center} 
\end{table}

The full automorphism group of the Doro graph, denoted by $G_{\Gamma_D}$, has order 16320. 
Table \ref{table-coeff-D} shows that the only possibility to construct self-orthogonal codes from $\Gamma_D$ by applying Corollary \ref{cor_drg} is a construction of 
codes over a field of characteristic 2 from the matrix $A_2$.   
The only subgroups of $G_{\Gamma_D}$ that act with all orbits of the same length are the trivial group and the groups of order 17 and 34, respectively. 
We constructed binary codes arising by applying Corollary \ref{cor_drg}, and in Table \ref{table-codesSO-D} we present the results obtained for the trivial group.

\begin{table}[H]
\begin{center} \begin{footnotesize} 
\begin{tabular}{|c|c|c|}
 \hline 
$H \leq G_{\Gamma_D}$ & $i$ & The code \\
\hline \hline
 $I$&2&$[68,8,32]_2$ $*$ \\
 \hline \hline
\end{tabular} \end{footnotesize}
\caption{\footnotesize Self-orthogonal codes obtained from the graph $\Gamma_D$} \label{table-codesSO-D}
\end{center} 
\end{table}

\begin{rem}
The code with parameters $[68,8,32]_2$ is a projective code with weights $32$ and $40$ which yields a strongly regular graph with parameters $(256,187,138,132)$, 
the complement of a strongly regular graph with parameters $(256,68,12,20)$.  
\end{rem}

\subsection{Hadamard graph on 48 vertices, $d=4$}

The Hadamard graph $\Gamma_H$ on 48 vertices has 48 vertices and the intersection array $\{12,11,6,1;1,6,11,12\}$. 
The coefficients of $A_i \cdot A_i$ for the distance-$i$ matrices $A_i$ of $\Gamma_H$ ($i=0, 1, \dots, 4$) are given in Table \ref{table-coeff-H}.

\begin{table}[H] 
\begin{center} \begin{scriptsize}
\begin{tabular}{|c|c|c|c|c|c|}
 \hline 
 & $A_0$ & $A_1$& $A_2$ & $A_3$ & $A_4$\\
 \hline \hline
$A_0 \cdot A_0$ & 1 & 0 & 0 & 0 &0\\
 \hline
$A_1 \cdot A_1$ & 12 & 0 & 6 & 0 &0 \\
 \hline
$A_2 \cdot A_2$ & 22 & 0 & 20 & 0 &22 \\
 \hline
 $A_3 \cdot A_3$ & 12 & 0 & 6 & 0 &0  \\
\hline
$A_4 \cdot A_4$ & 1 & 0 & 0 & 0 &0  \\
 \hline \hline
\end{tabular} \end{scriptsize}
\caption{\footnotesize The coefficients $p_{i,i}^k$ of adjacency matrices for the graph $\Gamma_H$} \label{table-coeff-H}
\end{center} 
\end{table}

The full automorphism group $G_{\Gamma_H}$ of $\Gamma_H$ has order 380160. 
Table \ref{table-coeff-H} shows that by Corollary \ref{cor_drg} one can construct binary codes from $A_1$, $A_2$ and $A_3$, and ternary codes form $A_1$ and $A_3$.
For the subgroups of $G_{\Gamma_H}$ acting with all orbits of the same length we construct the quotient matrices $M_i^{H}$, $i=1,2$, which span self-orthogonal codes. 
Since the distance graphs $ (\Gamma_{H})_1$ and $(\Gamma_{H})_3$ are isomorphic, we do not have to check quotient matrices of $(\Gamma_{H})_3$.
Table \ref{table-codesSO-H} presents the obtained results.

\begin{table}[H]
\begin{center} \begin{footnotesize} 
\begin{tabular}{|c|c|c||c|c|c|}
 \hline 
$H\leq G_{\Gamma_H}$ & $i$ & The code & $H\leq G_{\Gamma_H}$ & $i$ & The code \\
\hline \hline
 $I$&1&$[48,24,4]_2$  &$Z_3$&2&$[16,8,2]_2$ \\ \hline
$I$&2&$[48,24,2]_2$  &$Z_3$&1&$[16,8,4]_2$ \\ \hline
 $I$&1&$[48,34,4]_3$ &$Z_3$&1&$[16,4,6]_3$ \\ \hline
$Z_2$&1&$[24,2,12]_2$ &$E_4$&2&$[12,4,2]_2$ \\ \hline
$Z_2$&2&$[24,8,2]_2$ &$E_4$&1&$[12,5,2]_2$ \\ \hline
 $Z_2$&1&$[24,10,2]_3$ &$Z_4, E_4$&1&$[12,2,6]_3$ \\ \hline
$Z_2$&2&$[24,12,2]_2$ &$E_4$&1&$[12,3,6]_3$ \\ \hline
$Z_2$&1&$[24,12,4]_2$ &$Z_6, S_3$&1&$[8,2,4]_2$  \\ \hline
 $Z_2$&1&$[24,2,12]_3$ &$Z_6, S_3$&1&$[8,4,4]_2$ $*$ \\ \hline
$Z_2$&1&$[24,6,6]_3$ &$Z_6, S_3$&2&$[8,4,2]_2$ \\ \hline
 $Z_2$&1&$[24,7,12]_3$ $*$ &$Z_6, S_3$&1&$[8,2,6]_3$ $*$ \\
 \hline \hline
\end{tabular} \end{footnotesize}
\caption{\footnotesize Self-orthogonal codes obtained from the graph $\Gamma_H$} \label{table-codesSO-H}
\end{center} 
\end{table}

\begin{rem}
From the supports of all codewords of weight $4$ of the code with parameters $[24,12,4]_2$, in a similar way as described in \cite{sicily} we construct the triangular graph $T(12)$, 
\textit{i.e.} the strongly regular graph with parameters $(66,20,10,4)$, 
and from the codewords of weight $8$ a strongly regular graph with parameters $(495,238,109,119)$. Further, we obtain a distance-regular graph of diameter $4$, 
known as Johnson graph with $495$ vertices having the intersection array $\{32,21,12,5;1,4,9,16\}$. The code with parameters $[8,4,4]_2$ is the famous Hamming code.
\end{rem}

\subsection{Doubled Gewirtz graph, $d=5$}

The Doubled Gewirtz graph $\Gamma_{DG}$ has 112 vertices, the full automorphism group $G_{\Gamma_{DG}}$ of order 161280, and the intersection array $\{10,9,8,2,1;1,2,8,9,10\}$. 
The coefficients of $A_i \cdot A_i$ for the distance-$i$ matrices $A_i$ of $\Gamma_{DG}$ ($i=0, 1, \dots, 5$) are given in Table \ref{table-coeff-DG}.
The group $G_{\Gamma_{DG}}$ contains subgroups acting on the graph $\Gamma_{DG}$ with all orbits of the same length.
For the corresponding quotient matrices $M_i^{H}$ we obtain self-orthogonal codes presented in Table \ref{table-codesSO-DG}.

\begin{table}[H] 
\begin{center} \begin{scriptsize}
\begin{tabular}{|c|c|c|c|c|c|c|}
 \hline 
 & $A_0$ & $A_1$& $A_2$ & $A_3$ & $A_4$& $A_5$\\
 \hline \hline
$A_0 \cdot A_0$ & 1 & 0 & 0 & 0 &0&0\\
 \hline
$A_1 \cdot A_1$ & 10 & 0 & 2 & 0 &0 &0\\
 \hline
$A_2 \cdot A_2$ & 45 & 0 & 36 & 0 &36 &0 \\
 \hline
 $A_3 \cdot A_3$ & 45 & 0 & 36 & 0 &36 &0  \\
\hline
$A_4 \cdot A_4$ & 10 & 0 & 2 & 0 &0 &0  \\
\hline
$A_5 \cdot A_5$ & 1 & 0 & 0 & 0 &0 &0  \\
 \hline \hline
\end{tabular} \end{scriptsize}
\caption{\footnotesize The coefficients $p_{i,i}^k$ of adjacency matrices for the graph $\Gamma_{DG}$} \label{table-coeff-DG}
\end{center} 
\end{table}

\begin{table}[H]
\begin{center} \begin{footnotesize} 
\begin{tabular}{|c|c|c||c|c|c|}
 \hline 
$H\leq G_{\Gamma_{DG}}$ & $i$ & The code &$H\leq G_{\Gamma_{DG}}$ & $i$ & The code \\
\hline \hline
$I$&1,4&$[112,40,10]_2$ &$Z_2$&2,3&$[56,19,18]_3$ \\ \hline
$I$&2,3&$[112,38,18]_3$ &$E_4$&1,4&$[28,9,8]_2$ \\ \hline
$Z_2$&1,4&$[56,18,8]_2$ &$E_4$&2,3&$[28,9,12]_3$ $*$ \\ \hline
$Z_2$ &1,4&$[56,20,10]_2$ &$Z_7$&1,4&$[16,4,2]_2$  \\ \hline
$Z_2$&2,3&$[56,18,12]_3$ &$Z_7$&2,3&$[16,2,6]_3$ \\
 \hline \hline
\end{tabular} \end{footnotesize}
\caption{\footnotesize Self-orthogonal codes obtained from the graph $\Gamma_{DG}$} \label{table-codesSO-DG}
\end{center} 
\end{table}

\begin{rem}
From the supports of all codewords of weight $10$ of the code with parameters $[56,20,10]_2$ we construct the Sims-Gewirtz graph, the unique strongly regular graph with parameters $(56,10,0,2)$, 
and from the codewords of weight $14$, as explained in \cite{sicily}, we construct the unique strongly regular $(120,42,8,18)$ graph, known as the $L(3,4)$ graph, defined on Baer subplanes.
\end{rem}

\subsection{Incidence graph of $GH(3,3)$, $d=6$}

The incidence graph of $GH(3,3)$, which we will denote by $\Gamma_{I}$, has 728 vertices, the intersection array $\{4,3,3,3,3,3;1,1,1,1,1,4\}$,
and the full automorphism group $G_{\Gamma_{I}}$ of order 8491392. 
The coefficients of $A_i \cdot A_i$ for the distance-$i$ matrices $A_i$ of $\Gamma_{I}$ ($i=0, 1, \dots, 6$) are given in Table \ref{table-coeff-I}.
As in the previous examples, by applying Corollary \ref{cor_drg} we obtain self-orthogonal codes which are presented in Table \ref{table-codesSO-I}.

\begin{table}[H] 
\begin{center} \begin{scriptsize}
\begin{tabular}{|c|c|c|c|c|c|c|c|}
 \hline 
 & $A_0$ & $A_1$& $A_2$ & $A_3$ & $A_4$& $A_5$& $A_6$\\
 \hline \hline
$A_0 \cdot A_0$ & 1 & 0 & 0 & 0 &0&0&0\\
 \hline
$A_1 \cdot A_1$ & 4&0&1&0&0&0&0\\
 \hline
$A_2 \cdot A_2$ & 12 & 0 & 2 & 0 &1 &0 &0 \\
 \hline
 $A_3 \cdot A_3$ & 36 & 0 & 6 & 0 &2 &0 &4  \\
\hline
$A_4 \cdot A_4$ & 108 & 0 & 18 & 0 &33 &0 &32  \\
\hline
$A_5 \cdot A_5$ & 324 & 0 & 297 & 0 &288&0 &288 \\
\hline
$A_6 \cdot A_6$ & 243 & 0 & 162 & 0 &162&0 &162 \\
 \hline \hline
\end{tabular} \end{scriptsize}
\caption{\footnotesize The coefficients $p_{i,i}^k$ of adjacency matrices for the graph $\Gamma_I$} \label{table-coeff-I}
\end{center} 
\end{table}

\begin{table}[H]
\begin{center} \begin{footnotesize} 
\begin{tabular}{|c|c|c||c|c|c|}
 \hline 
$H\leq G_{\Gamma_I}$ & $i$ & The code &$H\leq G_{\Gamma_I}$ & $i$ & The code\\
\hline \hline
 $Z_7$&3&$[104,26,12]_2$ &$Z_{13}$&6&$[56,6,18]_3$ \\ \hline
 $Z_7$&5&$[104,26,18]_3$ &$D_{14}, Z_{14}$&3&$[52,13,12]_2$  \\ \hline
$Z_7$&6&$[104,6,36]_3$  &$D_{14}, Z_{14}$&5&$[52,13,18]_3$  \\ \hline
 $Z_{13}$&3&$[56,14,8]_2$ &$D_{14}, Z_{14}$&6&$[52,3,36]_3$ $*$ \\ \hline
$Z_{13}$&5&$[56,14,9]_3$ &&&\\
 \hline \hline
\end{tabular} \end{footnotesize}
\caption{\footnotesize Self-orthogonal codes obtained from the graph $\Gamma_{I}$} \label{table-codesSO-I}
\end{center} 
\end{table}

\subsection{Doubled Odd graph $D(O_4)$, $d=7$}

The Doubled Odd graph $D(O_4)$ has 70 vertices, the full automorphism group $G_{D(O_4)}$ of order 10080, and the intersection array $\{4,3,3,2,2,1,1;1,1,2,2,3,3,4\}$. 
The coefficients of the distance-$i$ matrices $A_i$ of $D(O_4)$ ($i=0, 1, \dots, 7$) are given in Table \ref{table-coeff-DO}.
By applying Corollary \ref{cor_drg}, in a similar way as in the previous examples, we obtain self-orthogonal codes and 
present the obtained results in Table \ref{table-codesSO-DO}.

\begin{table}[H] 
\begin{center} \begin{scriptsize}
\begin{tabular}{|c|c|c|c|c|c|c|c|c|}
 \hline 
 & $A_0$ & $A_1$& $A_2$ & $A_3$ & $A_4$& $A_5$& $A_6$& $A_7$\\
 \hline \hline
$A_0 \cdot A_0$ & 1 & 0 & 0 & 0 &0&0&0&0\\
 \hline
$A_1 \cdot A_1$ & 4&0&1&0&0&0&0&0\\
 \hline
$A_2 \cdot A_2$ & 12 & 0 & 5 & 0 &4 &0 &0 &0 \\
 \hline
 $A_3 \cdot A_3$ & 18 & 0 & 9 & 0 &9 &0 &9 &0  \\
 \hline
$A_4 \cdot A_4$ & 18 & 0 & 9 & 0 &9 &0 &9 &0  \\
 \hline
$A_5 \cdot A_5$ & 12 & 0 & 5 & 0 &4 &0 &0 &0 \\
 \hline
$A_6 \cdot A_6$ & 4&0&1&0&0&0&0&0 \\
 \hline
$A_7 \cdot A_7$ & 1&0&1&0&0&0&0&0 \\
 \hline \hline
\end{tabular} \end{scriptsize}
\caption{\footnotesize The coefficients $p_{i,i}^k$ of adjacency matrices for the graph $D(O_4)$} \label{table-coeff-DO}
\end{center} 
\end{table}

\begin{table}[H]
\begin{center} \begin{footnotesize} 
\begin{tabular}{|c|c|c|}
 \hline 
$H\leq G_{D(O_4)}$ & $i$ & The code \\
\hline \hline
 $I$&3,4&$[70,26,12]_3$ \\ \hline
 $Z_2$&3,4&$[35,13,12]_3$ \\ \hline
$Z_5$&3,4&$[14,2,6]_3$ \\ \hline
 $Z_7$&3,4&$[10,2,3]_3$ \\
 \hline \hline
\end{tabular} \end{footnotesize}
\caption{\footnotesize Self-orthogonal codes obtained from the graph $D(O_4)$} \label{table-codesSO-DO}
\end{center} 
\end{table}

\subsection{Foster graph, $d=8$}

The Foster graph $\Gamma_F$ has 90 vertices, the full automorphism group $G_{\Gamma_F}$ of order 4320, and the intersection array $ \{3,2,2,2,2,1,1,1; 1,1,1,1,2,2,2,3 \} $. 
The coefficients of the distance-$i$ matrices $A_i$ ($i=0, 1, \dots, 8$) are given in Table \ref{table-coeff-FG}.
In a similar way as in the previous examples, we obtain self-orthogonal codes by applying Corollary \ref{cor_drg}. 
The results are presented in Table \ref{table-codesGF-FG}.
 
\begin{table}[H] 
\begin{center} \begin{scriptsize}
\begin{tabular}{|c|c|c|c|c|c|c|c|c|c|}
 \hline 
 & $A_0$ & $A_1$& $A_2$ & $A_3$ & $A_4$ &$A_5$ & $A_6$& $A_7$ & $A_8$ \\
\hline \hline
$A_0 \cdot A_0$ & 1 & 0 & 0 & 0 & 0 & 0& 0& 0& 0 \\
 \hline
$A_1 \cdot A_1$ & 3 & 0 & 1 & 0 & 0 & 0& 0& 0& 0 \\
 \hline
$A_2 \cdot A_2$ & 6 & 0 & 1 & 0 & 1 & 0& 0& 0& 0 \\
 \hline
 $A_3 \cdot A_3$ & 12 & 0 & 2 & 0 & 3 & 0& 4& 0& 0 \\
 \hline
 $A_4 \cdot A_4$ & 24 & 0 & 12 & 0 & 12 & 0& 12& 0& 24 \\
 \hline
 $A_5 \cdot A_5$ & 24 & 0 & 12 & 0 & 12 & 0& 14& 0& 12 \\
 \hline
 $A_6 \cdot A_6$ & 12 & 0 & 2 & 0 & 4 & 0& 1& 0& 6 \\
 \hline
 $A_7 \cdot A_7$ & 6 & 0 & 2 & 0 & 0 & 0& 1& 0& 3 \\
 \hline
 $A_8 \cdot A_8$ & 2 & 0 & 0 & 0 & 0 & 0& 0& 0& 1 \\
 \hline \hline
\end{tabular} \end{scriptsize}
\caption{\footnotesize The coefficients $p_{i,i}^k$ of adjacency matrices for the Foster graph $\Gamma_F$} \label{table-coeff-FG}
\end{center} 
\end{table}

\begin{table}[H]
\begin{center} \begin{footnotesize} 
\begin{tabular}{|c|c|c||c|c|c|}
 \hline 
$H\leq G_{\Gamma_F}$ & $i$ & The code &$H\leq G_{\Gamma_F}$ &$i$ & The code\\
\hline \hline
$I$&5&$[90,12,20]_2$ &$Z_3$&4&$[30,8,8]_2$    \\ \hline
$I$&4& $[90,8,24]_2$  &$Z_5$&5&$[18,4,4]_2$    \\ \hline
$I$&4&$[90,30,3]_3$ &$Z_5$&4&$[18,6,3]_3$ \\ \hline
$Z_2$&5&$[45,6,20]_2$ &$S_3$&4&$[15,4,8]_2$ $*$\\ \hline
$Z_2$&4&$[45,4,24]_2$ &$D_{10}, Z_{10}$&5&$[9,2,4]_2$\\ \hline
$Z_2$&4&$[45,15,3]_3$ &$D_{10}, Z_{10}$&4&$[9,3,3]_3$\\
\hline \hline
\end{tabular} \end{footnotesize}
\caption{\footnotesize Self-orthogonal codes from the graph $\Gamma_F$} \label{table-codesGF-FG}
\end{center} 
\end{table}

\begin{rem}
The code with parameters $[45,6,20]_2$ is a projective code with weights $20$ and $24$ which yields a strongly regular graph with parameters $(64,45,32,30)$. 
The code with parameters $[30,8,8]_2$ is a projective code with weights $8$ and $16$ which yields the unique strongly regular graph with parameters $(256,30,14,2)$ 
and from the supports of codewords of weight $16$ in a way explained in \cite{sicily} we obtain the unique strongly regular graph with parameters $(225,28,13,2)$. 
From the supports of the codewords of weight $6$ of the code $[18,6,3]_3$ we obtain the unique strongly regular graph with parameters $(15,6,1,3)$, 
and from the supports of the codewords of weight $9$ we obtain a distance-regular graph with diameter $3$ having $20$ vertices and intersection array $\{9,4,1;1,4,9\}$. 
The code $[18,4,4]_2$ is not projective but it is two weight code with weights $4$ and $8$ which yields a strongly regular graph with parameters $(16,6,2,2)$. 
From the supports of codewords of weight $8$ of the same code we obtain the unique strongly regular graph with parameters $(9,4,1,2)$.
\end{rem}

\section{Self-orthogonal subspace codes from equitable partitions of association schemes}\label{SO-subspace}

As an analog of the definition of a self-orthogonal linear code we introduce the definition of a self-orthogonal subspace code as follows.

\begin{defi}
Let ${\cal P}_q(n)$ be the set of all subspaces of $\mathbb F_q^n$. The {\it dual} code of a subspace code $C_S \subseteq  {\cal P}_q(n)$ is the set $C_S^{\perp}$ of all vector spaces in 
${\cal P}_q(n)$ that are orthogonal to each vector space in $C_S$.
If $C_S \subseteq C_S^\perp$, then $C_S$ is called {\it self-orthogonal} subspace code. If $C_S = C_S^\perp$, then $C_S$ is called {\it self-dual}.
\end{defi}

Let $\mathcal{A}= \{A_0, A_1, ..., A_d \}$, be the set of $n \times n$ adjacency matrices of the association scheme $(X,\mathcal{R})$, and let $q=p^m$ be a prime power.
Let us consider the matrix algebra $\overline {\cal A}$ over $\mathbb F_q$ generated by matrices $A_{i_1}, A_{i_2}, ..., A_{i_t}$, $I=\{i_1, i_2, ..., i_t\} \subseteq \{0, 1, ..., d \}$, $t \geq 2$.
If $p|p_{x,y}^k$, for all $k \in \{0,1,...,d\}$ and all choices of $x, y \in I$, then the row spaces of the elements of $\overline {\cal A}$ are mutually orthogonal.
Therefore, the set of the row spaces of elements of $\overline {\cal A}$ forms a self-orthogonal subspace code $C_S \subseteq \mathbb F_q^n$.
Moreover, because of Theorem \ref{th_struc_const} the following theorem holds.

\begin{tm} \label{so_subcode}
Let $\Pi$ be an equitable partition of a $d$-class association scheme $(X,\mathcal{R})$ with $n$ cells of the same length $\frac{|X|}{n}$ and let $p$ be a prime number.
Further, let $I=\{i_1, i_2, ..., i_t\} \subseteq \{0, 1, ..., d \}$ and $p|p_{x,y}^k$, for all $k \in \{0,1,...,d\}$ and all $x, y \in I$.
Then the set of row spaces of elements of the matrix algebra generated by the matrices $M_i$, $i \in I$, forms a self-orthogonal subspace code $C_S \subseteq \mathbb F_q^n$,
where $q=p^m$ is a prime power.
\end{tm}

\subsection{Subspace codes from distance-regular graphs}

Let $\Gamma$ be a distance-regular graph with diameter $d$ and adjacency matrix $A$. Let $A_i$ denotes the distance-$i$ matrix of $\Gamma$, for $i=0, 1, \dots, d$. 
Let $G$ be an automorphism group of $\Gamma$. In the sequel, the quotient matrix of $A_i$, where $i=1, 2, \dots, d$,  with respect to the orbit partition induced by $G$ will be denoted by $M_i$.

Since a distance-regular graph induces an association scheme, the next corollary follows from Theorem \ref{so_subcode}.

\begin{kor} \label{cor_drg_sub}
Let $\Gamma$ be a distance-regular graph with diameter $d$, and  let an automorphism group $G$ acts on $\Gamma$ with $n$ orbits of the same length.
Further, let $I=\{i_1, i_2, ..., i_t\} \subseteq \{0, 1, ..., d \}$ and $p$ be a prime number such that $p|p_{x,y}^k$, for all $k \in \{0,1,...,d\}$ and all $x, y \in I$.
Then the set of row spaces of elements of the matrix algebra generated by the matrices $M_i$, $i \in I$, forms a self-orthogonal subspace code $C_S \subseteq \mathbb F_q^n$,
where $q=p^m$ is a prime power.
\end{kor}

\subsubsection{Examples of self-orthogonal subspace codes}

As an example, we apply Corollary \ref{cor_drg_sub} to a construction of self-orthogonal subspace codes from the Doubled Higman-Sims graph.
The Doubled Higman-Sims graph, which we denote by with $\Gamma_{dHS}$, is a distance-regular graph having 200 vertices, diameter $d=5$, and intersection array $\{22,21,16,6,1;1,6,16,21,22\}$. 
See \cite{BCN} for more information.
The coefficients of $A_i \cdot A_i$ for the distance-$i$ matrices $A_i$ of $\Gamma_{dHS}$ ($i=0, 1, \dots, 5$) are given in Table \ref{table-coeff-dHS}.

\begin{table}[H] 
\begin{center} \begin{scriptsize}
\begin{tabular}{|c|c|c|c|c|c|c|}
 \hline 
 & $A_0$ & $A_1$& $A_2$ & $A_3$ & $A_4$& $A_5$\\
 \hline \hline
$A_0 \cdot A_0$ & 1 & 0 & 0 & 0 &0&0\\
 \hline
$A_1 \cdot A_1$ & 22&0&6&0&0&0\\
 \hline
$A_2 \cdot A_2$ & 77 & 0 & 60 & 0 &56&0 \\
 \hline
 $A_3 \cdot A_3$ & 77 & 0 & 60 & 0 &56&0  \\
\hline
$A_4 \cdot A_4$ & 22&0&6&0&0&0  \\
\hline
$A_5 \cdot A_5$ & 1 & 0 & 0 & 0 &0&0  \\

 \hline \hline
\end{tabular} \end{scriptsize}
\caption{\footnotesize The coefficients $p_{i,i}^k$ of adjacency matrices for the graph $\Gamma_{dHS}$} \label{table-coeff-dHS}
\end{center} 
\end{table}

Table \ref{table-coeff-dHS} shows that we have to check only the case $I=\{1,4 \}$, \textit{i.e.} the remaining product $A_1\cdot A_4$. 
Since, $A_1\cdot A_4=6A_3+22A_5$, and $p=2$ divides all coefficients $p_{x,y}^k$ for $k \in \{0, 1, ...,5 \}$ and $x, y \in I$, 
the set of row spaces of the elements of the matrix algebra generated by matrices $M_i$, $i \in I$, is a self-orthogonal subspace code.

The full automorphism group $G_{\Gamma_{dHS}}$ of the  Doubled Higman-Sims graph has order 177408000. For the subgoups of $G_{\Gamma_{dHS}}$ acting on the graph 
$\Gamma_{dHS}$ with all orbits of the same length we construct quotient matrices $M_i$, $i \in I$. 
The self-orthogonal subspace codes obtained by applying Corollary \ref{cor_drg_sub} are presented in Table \ref{HSdrg}.

\begin{table}[H]
\begin{center} \begin{footnotesize}
\begin{tabular}{|c|c|}
 \hline
 $H\leq G_{\Gamma_{dHS}} $ & The code \\
\hline \hline
$E_{25}$ &$(8,8,1,\{0,1,2,3,4\})_2$ \\
\hline
$Z_5: Z_4$ &$(10,4,1,\{0,1,2\})_2$ \\
$ Z_5: Z_4$&$(10,5,1,\{0,1,2\})_2$ \\
\hline
$Z_{10}$&$(20,236,1,\{0,1,2,3,4,5,6\})_2$ \\
$Z_{10}$&$(20,107,1,\{0,1,2,3,4,5,6\})_2$ \\
$Z_{10}$&$(20,67,1,\{0,1,2,3,4\})_2$ \\
$Z_{10}$&$(20,3,2,\{0,1,2\})_2$ \\
\hline
$Z_4\times  Z_2$&$(25,366,1,\{0,1,2,3,4,5\})_2*$ \\
$D_8$&$(25,83,1,\{0,1,2,3,4,5\})_2$ \\
$E_8$&$(25,67,1,\{0,1,2,3,4\})_2$ \\
$E_8$&$(25,17,1,\{0,1,2,3,4\})_2$ \\
$D_8$&$(25,6,1,\{0,1,2,3,4\})_2$ \\
\hline
$Z_5$&$(40,27,1,\{0,1,2,3,4\})_2$ \\
\hline
$I$ &$(200,3,22,\{0,22,44\})_2$\\
\hline \hline
\end{tabular} \end{footnotesize}
\caption{\footnotesize Self-orthogonal subspace codes from the graph $\Gamma_{dHS}$} \label{HSdrg}
\end{center}
\end{table}

\begin{rem}
One of the codewords (\textit{i.e.} subspaces) of the subspace code labeled with $*$ is an optimal linear code with parameters $[25,4,12]_2$, 
a two weight code with weights $12$ and $16$. It is not projective but, it yields the unique strongly regular graph with parameters $(16,10,6,6)$, 
the complement of the strongly regular graph with parameters $(16,5,0,2)$.
From the supports of the codewords of weight $12$ one obtains the Peterson graph, the unique strongly regular graph with parameters $(10,3,0,1)$.
\end{rem}

\section{Conclusion} \label{con}

In this paper we gave a method of constructing self-orthogonal codes from equitable partitions of association schemes and, as an illustration of the method, 
constructed self-orthogonal codes from some distance-regular graphs.
Some of the obtained codes are optimal.
We also introduced self-orthogonal subspace codes and showed that under some conditions equitable partitions of association schemes yield such subspace codes.
We gave examples of self-orthogonal subspace codes obtained from distance-regular graph. We believe that using the described methods one can obtain further interesting self-orthogonal codes
and self-orthogonal subspace codes using association schemes, especially those arising from distance-regular graphs.


\begin{thebibliography}{30}

\bibitem{magma}
Bosma, W., Cannon, J.: Handbook of Magma Functions, Department of Mathematics, University of Sydney, 1994. (http://magma.maths.usyd.edu.au/magma)

\bibitem{BCN}
Brouwer, A. E., Cohen, A. M., Neumaier, A.: Distance-Regular Graphs. Springer-Verlag, Berlin, (1989)

\bibitem{SRG-codes}
Crnkovi\'c, D., Maksimovi\'c, M., Rodrigues, B. G., Rukavina, S.: Self-orthogonal codes from the strongly regular graphs on up to 40 vertices, Adv. Math. Commun. 10 (2016), 555--582

\bibitem{sicily}
Crnkovi\'c, D., Mikuli\'c Crnkovi\'c, V., \v Svob, A.: Transitive combinatorial structures invariant under some subgroups of $S(6,2)$ and related codes, Atti Accad. Peloritana Pericolanti Cl. Sci. Fis. Mat. Natur. 96, No. S2, A6 (2018), 15 pages.

\bibitem{2-des-codes}
Crnkovi\'c, D., Rodrigues, B. G., Rukavina, S., Sim\v ci\'c, L.: Self-orthogonal codes from orbit matrices of 2-designs, Adv. Math. Commun. 7 (2013), 161--174

\bibitem{GoM}
Godsil, C. D., Martin, W. J.: Quotients of Association Schemes, J. Combin. Theory Ser. A 69 (1995), 185--199

\bibitem{codetables}
Grassl, M.: Bounds on the minimum distance of linear codes and quantum codes, (http://www.codetables.de.)

\bibitem{network-book}
Greferath, M., Pav\v cevi\'c, M. O., Silberstein, N., Vázquez-Castro, M. \'A. (eds.)
Network coding and subspace designs. Signals and Communication Technology. Springer, Cham, (2018). 

\bibitem{Hanaki}
Hanaki, A., Elementary functions for association schemes on GAP, October 2013. (http://math.shinshu-u.ac.jp/$\sim$hanaki/as/gap/association\_scheme.pdf)

\bibitem{har_ton} Harada, M., Tonchev, V. D.: Self-orthogonal codes from symmetric designs with fixed-point-free automorphisms, 
Discrete Math. 264 (1-3) (2003), 81--90

\bibitem{table_sc}
Heinlein, D., Kiermaier, M., Kurz, S., Wassermann, A.: Tables of subspace codes, arXiv: 1601.02864v2, 2017.

\bibitem{heinlein-dcc}
Heinlein, D., Honold, T., Kiermaier, M., Kurz S., Wassermann A.: Classifying optimal binary subspace codes of length 8, constant dimension 4 and minimum distance 6, Des. Codes Cryptogr. (2019), 87--375. https://doi.org/10.1007/s10623-018-0544-8

\bibitem{FEC} 
Huffman, W. C., Pless, V.: Fundamentals of Error-Correcting Codes, Cambridge University Press, 2003.

\bibitem{network-coding}
K\"{o}tter, R., Kschischang, F.: Coding for errors and erasures in random network coding, 
IEEE Trans. Inform. Theory 54 (2008), 3579–3591.

\bibitem{BMTh}
Martin, W. J.: Completely Regular Subsets, Ph. D: Thesis, University of Waterloo, April 1992.

\bibitem{RainsSloane}
Nebe, G., Rains,E. M., Sloane, N. J. A.: Self-Dual Codes and Invariant Theory,  Algorithms and Computation in Mathematics, Vol. 17, Springer-Verlag, 2006.


\bibitem{sage}
SageMath, the Sage Mathematics Software System (Version 8.2), The Sage Developers, 2018 (http://www.sagemath.org)

\bibitem{GAP2016}
The GAP Group, GAP -- Groups: Algorithms, and Programming, Version 4.8.4; 2016. (http://www.gap-system.org) 

\end{thebibliography}
\end{document}